\documentclass[11pt,a4paper]{article}

\usepackage{epsf,epsfig,amsfonts,amsgen,amsmath,amstext,amsbsy,amsopn,amsthm
}
\usepackage{amsmath,times,mathptmx}
\usepackage{amsfonts,amsthm,amssymb}
\usepackage{amsfonts}
\usepackage{graphics}
\usepackage{latexsym,bm}
\usepackage{amsfonts,amsthm,amssymb,bbding}
\usepackage{indentfirst}
\usepackage{graphicx}
\usepackage{color}
\usepackage[colorlinks=true,anchorcolor=blue,filecolor=blue,linkcolor=blue,urlcolor=blue,citecolor=blue]{hyperref}
\usepackage{float}
\usepackage{tikz}
\newtheorem{thm}{Theorem}

\newtheorem{lemma}{Lemma}
\newtheorem{false statement}{False statement}

\theoremstyle{definition}

\newtheorem{claim}{Claim}

\newtheorem{case}{Case}

\newtheorem{problem}{Problem}

\begin{document}
\allowdisplaybreaks[4]
\title{Essential connectivity and spectral radius of graphs
\footnote{Supported by NSFC (Nos. 12361071 and 11901498).}}
\author{{Wenxiu Ding, Dan Li}\thanks{Corresponding author. E-mail: ldxjedu@163.com.}{, Yu Wang}{, Jixiang Meng}\\
{\footnotesize College of Mathematics and System Science, Xinjiang University, Urumqi 830046, China}}
\date{}

\maketitle {\flushleft\large\bf Abstract:}
A graph is trivial if it contains one vertex and no edges. The essential connectivity $\kappa^{\prime}$ of $G$ is defined to be the minimum number of vertices of $G$ whose removal produces a disconnected graph with at least two non-trivial components. Let $\mathcal{A}_n^{\kappa',\delta}$ be the set of graphs of order $n$ with minimum degree $\delta$ and essential connectivity $\kappa'$. In this paper, we determine the graphs attaining the maximum spectral radii among all graphs in $\mathcal{A}_n^{\kappa',\delta}$ and characterize the corresponding extremal graphs.
In addition, we also determine the digraphs which achieve the maximum spectral radii among all strongly connected digraphs with given essential connectivity and give the exact values of the spectral radii of these digraphs.

\vspace{0.1cm}
\begin{flushleft}
\end{flushleft}

\section{Introduction}
Throughout this paper we only consider simple connected graphs. Let $G=(V(G),E(G))$ be a graph with vertex set $V(G)$ and edge set $E(G)$, where $|V(G)|=n$. The adjacency matrix $A(G)$ of $G$ is an $n\times n$ matrix whose $(i,j)$ entry is 1 if $v_{i}$ is adjacent to $v_{j}$ in $G$, and 0 otherwise. It is obvious that $A(G)$ is a real symmetric matrix. Thus its eigenvalues are all real numbers. Let $\lambda_i(G)$ denote the $i$th largest eigenvalue of $A(G)$ for $i=1,2,\ldots,n$. Specifically, the largest eigenvalue of $A(G)$ is called the spectral radius of $G$, denoted by $\rho(G)$. It follows from the Perron-Frobenius Theorem that $\rho(G)$ is a simple eigenvalue of $A(G)$ and there exists a unique positive unit vector $\bm{x}$ called Perron vector such that $A(G)\bm{x}=\rho(G)\bm{x}$. Let $N_{G}(u)=\{v\in V(G):uv\in E(G)\}$ be the set of neighbors of $u$ in $G$ and $N_G[u]=N_G(u)\cup \{u\}$. The degree of the vertex $u$ is $d_G(u)=|N_{G}(u)|$. A graph $G$ is called $d$-regular if every vertex has the same degree equal to $d$. Denote by $\delta(G)$ (or for short $\delta$) and $\Delta(G)$ the minimum and maximum degrees of the vertices of $G$. Moreover, we know $\delta(G)\leq \rho(G)\leq \Delta(G)$, and the equality holds in either of these inequalities if and only if $G$ is regular. Denote by $D(G)$ the diagonal matrix of vertex degree of $G$. The Laplacian matrix of $G$ is defined as $L(G)=D(G)-A(G)$. Let $\mu_i(G)$ be the $i$th smallest eigenvalue of $L(G)$. Evidently, $L(G)$ is a positive semidefinite matrix and $\mu_n(G)\geq \cdots\geq\mu_2(G)\geq\mu_1(G)=0$. For a $d$-regular graph on $n$ vertices, $\lambda_i=d-\mu_i$ for $1\leq i\leq n$.

Connectivity is a fundamental concept in graph theory.
The connectivity (respectively, edge-connectivity) of a graph $G$ is the minimum number of vertices (respectively, edges) whose deletion results in a disconnected graph (or reduces it to a single vertex in case of connectivity). We use $\kappa(G)$ and $\lambda(G)$ to represent the connectivity and edge-connectivity of a graph $G$, respectively.
The study of connectivity has many applications in measuring reliability and fault-tolerance networks \cite{M.J. Ma}. A network can be modeled as a simple graph $G$.
In general, the larger $\kappa(G)$ or $\lambda(G)$ is, the more reliable the network is.
It is well known that $\kappa(G)\leq \lambda(G)\leq \delta(G)$.
A graph $G$ is called maximally edge connected or $\lambda$-optimal if $\lambda(G)=\delta(G)$ and maximally vertex connected if $\kappa(G)=\delta(G)$. However, $\kappa(G)$ or $\lambda(G)$ always underestimates the resilience of large networks. To overcome such shortcoming, Harary \cite{F. Harary} introduced the concept of conditional connectivity by imposing some conditions on the components of $G-S$, where $S$ is a subset of edges or vertices.
A graph is trivial if it contains one vertex and no edges. A vertex (edge-) cut $S$ of $G$ is essential if $G-S$ has at least two non-trivial components. For an integer $k>0$, a graph $G$ is essentially $k$-(edge)-connected if $G$ does not have an essential vertex (edge-) cut $S$ with $|S|<k$. In particular, the essential connectivity of $G$, denoted by $\kappa^{\prime}(G)$ or $\kappa^{\prime}$, is the minimum cardinality over all essential vertex cut of $G$. Essential connectivity, as a conditional connectivity, has improved the accuracy of reliability and fault tolerance analysis of networks, and has been widely used in the analysis of various networks.

This paper focuses on the study from spectral perspectives. The connections between the connectivity of graphs and eigenvalues of associated matrices have been well studied in past several decades.
Fiedler \cite{M. Fiedler} proved that $\mu_2(G)\leq \kappa(G)$ for a non-complete simple graph $G$. This seminal result provided researchers with another parameter that quantitatively measures the connectivity of a graph. Hence, $\mu_2(G)$ is known as the ``algebraic connectivity'' of $G$. A lot of research in graph theory over the last forty years was stimulated by Fiedler's work. This results have been improved and refined by various authors in \cite{A. Abiad, S.L. Chandran, S.M. Cioaba}. In particular, the spectral conditions for the connectivity $\kappa(G)$ of graphs have been extensively investigated. Berman and Zhang \cite{A. Berman} studied the spectral radius of graphs with $n$ vertices and $k$ cut vertices and described the graph that has the maximal spectral radius in this class. The union $G_1\cup G_2$ is defined to be $G_1\cup G_2=(V_1\cup V_2,~E_1\cup E_2)$. The join $G_1\vee G_2$ is obtained from $G_1\cup G_2$ by adding all the edges joining a vertex of $G_1$ to a vertex of $G_2$. $K_k\vee(K_1\cup K_{n-k-1})$ is shown to be the graph with the maximal spectral radius among all graphs of order $n$ with connectivity $\kappa(G)\leq k$ in \cite{J. Li, M.L. Ye}. Lu and Lin \cite{H.L. Lu} confirmed that $K_k\vee(K_{\delta-k+1}\cup K_{n-\delta-1})$ is the graph with the maximum spectral radius among all graphs of order $n$ with minimum degree $\delta(G)\geq k$ and connectivity $\kappa(G)\leq k$. The $l$-connectivity $\kappa_l$ of a graph $G$ is the minimum number of vertices of $G$ whose removal produces a disconnected graph with at least $l$ components or a graph with fewer than $l$ vertices. Fan, Gu and Lin \cite{D.D. Fan3} characterized the graphs which have the maximal spectral radius among all graphs of order $n$ with given minimum degree $\delta$ and $l$-connectivity $\kappa_l$.
The spectral conditions for the edge-connectivity $\lambda(G)$ of graphs have also been well investigated by many researchers, see \cite{D.D. Fan2, W.J. Ning}.

Since the essential connectivity is considered as an extension of the connectivity, which can better evaluate the fault tolerability of networks. Motivated by the above spectral results, we are interested in studying the spectral condition for essential connectivity. Thus, we pose the following problem.
\begin{problem}\label{p1}
Which graphs attain the maximum spectral radii among all connected graphs of order $n$ with fixed minimum degree $\delta$ and essential connectivity $\kappa^{\prime}$?
\end{problem}

Let $\mathcal{A}_n^{\kappa',\delta}$ be the set of graphs of order $n$ with minimum degree $\delta$ and essential connectivity $\kappa'$. Denote by $K_n$ the complete graph with order $n$.
Note that there is no bound between $\kappa'$ and $\delta$. For $\kappa'>\delta-1$ and an isolated vertex $u$, let $G_n^{\kappa',\delta}$ be the graph obtained from $\{u\}\cup (K_{\kappa'}\vee(K_1\cup K_{n-\kappa'-2}))$ by adding $\delta-1$ edges between $u$ and $K_{\kappa'}$, and then adding an edge between $u$ and $K_1$. If $\kappa'\leq\delta-1$, let $G_n^{\kappa',\delta}=K_{\kappa'}\vee(K_{n-\delta-1}\cup K_{\delta-\kappa'+1})$. In this paper, we mainly give the answer to Problem \ref{p1}. Then we have the following result.

\vspace*{2mm}
\begin{thm} \label{thm1}
Let $G$ be a connected graph of order $n\geq\kappa'+4$ with minimum degree $\delta$ and essential connectivity $\kappa'$. Then
$$\rho(G)\leq\rho(G_n^{\kappa',\delta}),$$
with equality if and only if $G\cong G_n^{\kappa',\delta}$.
\end{thm}

We also consider an analogue for digraphs.
For a digraph $D=(V(D), E(D))$, where $V(D)$ and $E(D)$ are the vertex set and arc set of $D$, respectively. Two vertices are called adjacent if they are connected by an arc.
If $e = uv\in E(D)$, then $u$ is the initial vertex of $e$ and $v$ is the terminal vertex.
Let $N_D^-(v)=\{u\in V(D)|(u,v)\in E(D)\}$ and $N_D^+(v)=\{u\in V(D)|(v,u)\in E(D)\}$ denote the in-neighbors and out-neighbors of $v$, respectively. Let $d_v^-=|N_D^-(v)|$ denote the indegree of the vertex $v$ and $d_v^+=|N_D^+(v)|$ denote the outdegree of the vertex $v$ in $D$.
We call $D$ strongly connected if for every pair $u, v \in V(D)$, there exists a directed path from $u$ to $v$ and a directed path from  $v$ to $u$. Let $\vec{K}_n$ denote the complete digraph on $n$ vertices in which there are
arcs $v_iv_j\in E(\vec{K}_n)$ and $v_jv_i\in E(\vec{K}_n)$ for every pair of distinct vertices $v_i,v_j\in V(\vec{K}_n)$. The essential connectivity of a strongly connected digraph is the minimum number of vertices whose removal produces a digraph with at least two strongly connected non-trivial components.\par

Given a strongly connected digraph $D$, let $A(D)=(a_{ij})$ denote the adjacency matrix of $D$, where $a_{ij}=1$ if $v_iv_j\in E(D)$ and $a_{ij} = 0$ otherwise, and let $\rho(D)$ denote its spectral radius, the largest eigenvalue of $A(D)$.
Over the past several decades, the spectral radius of digraphs with some given parameters have been well investigated. Brualdi \cite{R.A. Brualdi} surveyed some results on the spectra of digraphs.
A vertex set $A\subseteq V(D)$ is acyclic if the induced subdigraph $D[A]$ is acyclic.
A partition of $V(D)$ into $k$ acyclic sets is called a $k$-coloring of $D$. The minimum integer $k$ for which there exists a $k$-coloring of $D$ is the dichromatic number $\chi(D)$ of the digraph $D$.
Mohar \cite{B. Mohar} gave a lower bound of spectral
radius of digraphs with given dichromatic number. Lin and Shu \cite{H.Q. Lin1} characterized the digraph which has the maximal spectral radius among all strongly connected digraphs with given order and dichromatic number. Drury and Lin \cite{S. Drury} determined the digraphs which have the minimum and second minimum spectral radius among all strongly connected digraphs with given order and dichromatic number.
$F$ is called a clique of $D$ if $D[F]$ is a complete subdigraph of $D$. The clique number of a digraph $D$ is the number of vertices in the largest clique of $D$.
The girth of $D$ is the length of the shortest directed cycle of $D$.
Lin, Shu, Wu and Yu \cite{H.Q. Lin2} characterized the extremal digraphs with minimum spectral radius among all digraphs with given clique number and girth, and the extremal digraphs with maximum spectral radius among all digraphs with given vertex connectivity.
Our result is closely related to and motivated by the aforementioned results obtained in \cite{H.Q. Lin2}. We are interested in investigating the problem of maximizing the spectral radius of digraphs with given essential connectivity.

\begin{figure}[H]
\begin{center}
\includegraphics[width=5.8cm, height=3cm]{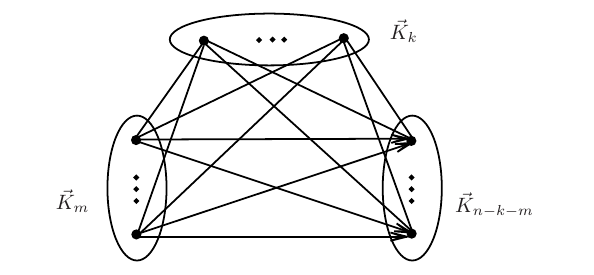}
\caption{$\vec{G}_n^{k,m}$}
\label{fig1}
\end{center}
\end{figure}

Let $\mathcal{D}_{n,k}$ be the set of strongly connected digraphs with order $n \geq \kappa'+4$ and essential connectivity $\kappa'(D)=k$.
Denote by $D_1 \nabla D_2$ the digraph obtained from two disjoint digraphs $D_1$ and $D_2$ with vertex set  $V(D_1) \cup V(D_2)$ and arc set $E=E(D_1) \cup E(D_2) \cup\{u v, v u \mid u \in V(D_1),~v \in V(D_2)\}$. Let $\vec {G}_n^{k,m}=\vec{K}_k \nabla(\vec{K}_m \cup \vec{K}_{n-k-m}) \cup F$, where $F=\{v u \mid v \in V(\vec{K}_m),~u \in V(\vec{K}_{n-k-m})\}$, as shown in Figure \ref{fig1}. Let $\vec{\mathcal{G}}(n,k)=\{\vec {G}_n^{k,m}\mid 2\leq m\leq n-k-2\}$. Obviously, $\vec{\mathcal{G}}(n,k)\subseteq \mathcal{D}_{n,k}$. Then we have the following result.

\begin{thm}\label{thm2}
 The digraphs $\vec {G}_n^{k,2}$ and $\vec {G}_n^{k,n-k-2}$ maximize the spectral radius among all digraphs in  $\mathcal{D}_{n, k}$.
\end{thm}

The current paper is organized as follows. In Section \ref{2}, we determine the graphs attaining the maximum spectral radii among all graphs with given minimum degree and essential connectivity. In addition, the corresponding extremal graphs are characterized. In Section \ref{3}, we characterize the digraphs which maximize the spectral radius with given essential connectivity, and we also give the exact values
of the spectral radii of these digraphs.

\section{Proof of Theorem \ref{thm1}} \label{2}
\begin{lemma} \cite{W.J. Ning} \label{lem1}
Let $u,v$ be two distinct vertices of a connected graph $G$, and let $\bm{x}$ be the Perron vector of $A(G)$.
\begin{description}
   \item[(1)] If $ N_G(v)\setminus \{u\}\subset N_G(u)\setminus \{v\}$, then $x_u>x_v$;
   \item[(2)] If $N_G(v)\subseteq N_G[u]$ and $N_G(u)\subseteq N_G[v]$, then $x_u=x_v$.
 \end{description}
\end{lemma}

\begin{lemma} \cite{Q. Li} \label{lem2}
Let $G$ be a connected graph and $G^{\prime}$ be a proper subgraph of $G$. Then $\rho(G^{\prime})<\rho(G)$.
\end{lemma}

\begin{lemma} \cite{M.L. Ye} \label{lem3}
Let $G$ be a connected graph and $\rho(G)$ be the spectral radius of $A(G)$. Let $u,v$ be two vertices of $G$. Suppose that $v_1,v_2, \ldots,v_s \in N_G(v) \setminus N_G(u)$ with $1\leq s\leq d_G(v)$, and $G^{\ast}$ is the graph obtained from $G$ by deleting the edges $v v_i$ and adding the edges $u v_i$ for $1 \leq i \leq s$. Let $\bm{x}$ be the Perron vector of $A(G)$. If $x(u) \geq x(v)$, then $\rho(G)<\rho(G^{\ast})$.
\end{lemma}

\begin{lemma}\label{lem4}\cite{D.D. Fan1}
Let $n=\sum_{i=1}^{t}n_i+s$. If $n_1\geq n_2\geq \cdots \geq n_t\geq p$ and $n_1<n-s-p(t-1)$, then
$$\rho(K_s\vee(K_{n_1}\cup K_{n_2} \cup\cdots \cup K_{n_t}))<\rho(K_s\vee(K_{n-s-p(t-1)}\cup (t-1)K_p)).$$
\end{lemma}

\noindent
{\bf{Proof of Theorem \ref{thm1}}.} Suppose that $G$ is a connected graph attaining the maximum spectral radius in $\mathcal{A}_n^{\kappa',\delta}$. By the definition of essential connectivity, there exists some nonempty subset $S\subseteq V(G)$ with $|S|=\kappa'$ such that $G-S$ contains at least two non-trivial components. Let $B_1,B_2,\ldots,B_q$ be the components of $G-S$ and at least two of them be non-trivial, where $q\geq 2$, and let $|V(B_i)|=n_i$ for $1\leq i\leq q$. Now, we divide the proof into the following two cases.
\vspace*{2mm}
\begin{case}$\kappa^{\prime}>\delta-1$. For any vertex $v\in V(G)$ and any subset $S\subseteq V(G)$, let $N_S(v)=N_G(v)\cap S$ and $d_S(v)=|N_G(v)\cap S|$. Choose a vertex $u\in V(G)$ such that $d_G(u)=\delta$ and suppose that $d_S(u)=t$.
\end{case}

\begin{claim}\label{claim1} $u\notin S$.
\end{claim}
Otherwise, $u\in S$. If $t=\delta$, then $N_G(u)\subseteq S$. This means that $\kappa'\geq\delta+1$. Let $S'=S-\{u\}$. Then $|S'|=\kappa'-1$, $G-S'$ contains $q+1$ components $u,B_1,B_2,\ldots,B_q$ and at least two of them are non-trivial components, contradicting the essential connectivity $\kappa'$ of $G$. If $t=\delta-1$, then $|N_G(u)\cap (V(G)\setminus S)|=1$. Without loss of generality, assume that $d_{B_1}(u)=1$. Let $S''=S-\{u\}$. Then $|S''|=\kappa^{\prime}-1$, $G-S''$ contains $q$ components $B_1\cup \{u\},B_2,\ldots,B_q$ and at least two of them are non-trivial components, also contradicting the essential connectivity $\kappa'$ of $G$. Thus we may assume $t\leq\delta-2$.
By the maximality of $\rho(G)$, we can deduce that $G-u\cong K_{\kappa^{\prime}-1}\vee(K_{n_1}\cup K_{n_2}\cup\cdots\cup K_{n_q})$. Note that $d_S(u)=t$ and $|N_G(u)\cap(V(G)\setminus S)|\geq2$. Let $N_S(u)=\{w_1,w_2,\ldots,w_t\}$ and $S\setminus N_S[u]=\{w_{t+1},w_{t+2},\ldots,w_{\kappa^{\prime}-1}\}$, and let $P_i=N_G(u)\cap V(K_{n_i})$ and $|P_i|=p_i$ for $1\leq i\leq q$. Then $\sum_{i=1}^{q}p_i=\delta-t$ and $\sum_{i=2}^{q}n_i\geq q$ because there is at least one non-trivial component of $B_2,B_3,\ldots,B_q$. Let $\bm{x}$ be the Perron vector of $A(G)$. By symmetry, $x(v)=x_i$ for $v\in P_i$ and $x(v)=x_i'$ for $v\in V(K_{n_i})\setminus P_i$, where $1\leq i\leq q$. Without loss of generality, assume that $x_1=\text{max}\{x_i|~1\leq i\leq q\}$. Observe that $x(w_1)=x(w_i)$ for $2\leq i\leq t$ and $x(w_{t+1})=x(w_i)$ for $t+2\leq i\leq \kappa^{\prime}-1$. Since $N_G(w_{t+1})\setminus \{w_1\}\subset N_G(w_1)\setminus \{w_{t+1}\}$, it follows that $x(w_1)>x(w_{t+1})$ by Lemma \ref{lem1}. By $A(G)\bm{x}=\rho(G)\bm{x}$, we have
\begin{equation}\label{equ1}
\begin{split}
\rho(G)x(u)&=tx(w_1)+\sum_{i=1}^{q} p_i x_i,
\end{split}
\end{equation}
\begin{equation}\label{equ2}
\begin{split}
\rho(G)x_i&=tx(w_1)+(\kappa^{\prime}-t-1) x(w_{t+1})+(p_i-1)x_i+(n_i-p_i)x_i^{\prime}+x(u),
\end{split}
\end{equation}
\begin{equation}\label{equ3}
\begin{split}
\rho(G)x_i^{\prime}&=tx(w_1)+(\kappa^{\prime}-t-1) x(w_{t+1})+p_ix_i+(n_i-p_i-1)x_i^{\prime},
\end{split}
\end{equation}
\begin{equation}\label{equ4}
\begin{split}
\rho(G)x(w_{t+1})&=tx(w_1)+(\kappa^{\prime}-t-2)x(w_{t+1})+\sum_{i=1}^{q} p_i x_i+\sum_{i=1}^{q}(n_i-p_i) x_i^{\prime},
\end{split}
\end{equation}
where $1\leq i \leq q$. From $(\ref{equ1})$, $(\ref{equ2})$ and $(\ref{equ3})$, we get
\begin{equation*}
\begin{aligned}
& \rho(G)\left(\sum_{i=2}^{q}p_ix_i+\sum_{i=2}^{q}(n_i-p_i)x_i^{\prime}-x(u)\right) \\
=&\left(\sum_{i=2}^{q}n_i-1\right)tx(w_1)+\sum_{i=2}^{q}n_i(\kappa^{\prime}-t-1)x(w_{t+1})+\sum_{i=2}^{q}(n_i-1)p_i x_i+\sum_{i=2}^{q}(n_i-1)(n_i-p_i) x_i^{\prime}\\
&+\sum_{i=2}^{q}p_ix(u)-\sum_{i=1}^{q}p_ix_i \\
> &(q-1)tx(w_1)+q(\kappa^{\prime}-t-1)x(w_{t+1})-\sum_{i=1}^{q}p_ix_i ~~(\text{since}~\sum_{i=2}^{q}n_i\geq q ~\text{and}~0 \leq p_i\leq n_i~\text{for}~2\leq i\leq q) \\
\geq &tx(w_1)+2(\kappa^{\prime}-t-1)x(w_{t+1})-(\delta-t)x_1~(\text{since}~q\geq 2,\sum_{j=1}^{q}p_j=\delta-t~\text{and}~x_1\geq x_i~\text{for}~2\leq i\leq q) \\
\geq&tx(w_{t+1})+2(\kappa^{\prime}-t-1)x(w_{t+1})-(\delta-t)x_1~~(\text{since}~x(w_1)>x(w_{t+1})~
\text{and}~t\geq 0)\\
=&(\kappa'-t)x(w_{t+1})+(\kappa'-2)x(w_{t+1})-(\delta-t)x_1\\
\geq &(\delta-t)(x(w_{t+1})-x_1)~~(\text{since}~\kappa^{\prime}\geq \delta\geq t+2~\text{and}~t\geq 0) .
\end{aligned}
\end{equation*}
Combining this with (\ref{equ2}) and (\ref{equ4}), we have
$$(\rho(G)+1)(x(w_{t+1})-x_{1})=\sum_{i=2}^{q}p_ix_i+\sum_{i=2}^{q}(n_i-p_i)x_i^{\prime}-x(u) >\frac{(\delta-t)(x(w_{t+1})-x_1)}{\rho(G)},$$
from which we have
$$\frac{(\rho^{2}(G)+\rho(G)-\delta+t)(x(w_{t+1})-x_1)}{\rho(G)} > 0.$$
Note that $G-S$ contains at least two non-trivial components. Without loss of generality, we may consider that
$B_2$ is one of the non-trivial components. And $K_{\kappa^{\prime}-1+n_2}$ is a proper subgraph of $G$. Thus, $\rho(G)>\rho(K_{\kappa^{\prime}-1+n_2})=\kappa^{\prime}+n_2-2 \geq \delta$ due to $n_2\geq2$ and $\kappa^{\prime}\geq \delta$. Combining this with $t\geq0$, we can deduce that $\rho^{2}(G)+\rho(G)-\delta+t>0$. This demonstrates that $x(w_{t+1})>x_1$. Recall that $\sum_{i=1}^{q}p_i=\delta-t\geq 2$ and $x_1=\text{max}\{x_i|~1\leq i\leq q\}$, then $p_1\geq 1$. Moreover, we can infer that $B_1$ is a non-trivial component of $G$. If not, we have $p_1=n_1=1$ and $p_2=0$ by the maximality of $x_1$. Then from (\ref{equ1}), (\ref{equ2}) and (\ref{equ3}), we have
\begin{equation*}
\begin{aligned}
&(\rho(G)-1)(x_2'-x_1)\\
=&(n_2-2)x_2'-x(u)+x_1\\
\geq& \frac{1}{\rho(G)+1}[(\rho(G)+1)(x_1-x(u))]~~(\text{since}~n_2\geq 2)\\
=&\frac{(\kappa'-t-1)x(w_{t+1})+x_1-\sum_{i=1}^{q}p_ix_i}{\rho(G)+1}\\
\geq& \frac{(\delta-t-1)(x(w_{t+1})-x_1)}{\rho(G)+1}~~(\text{since}~\kappa'\geq\delta,~\sum_{j=1}^{q}p_j=\delta-t~\text{and}~x_1\geq x_i~\text{for}~1\leq i\leq q)\\
>&0~~(\text{since}~t\leq\delta-2~\text{and}~x(w_{t+1})>x_1).
\end{aligned}
\end{equation*}
Since $\rho(G)>\delta\geq 1$, it follows that $x_2'>x_1$. Denote by $V(B_1)=\{v_1\}$ and we take a vertex $u_1\in V(K_{n_2})$. Let $G^*=G-\{uv_1\}+\{uu_1\}$, then $\rho(G^*)>\rho(G)$ by Lemma \ref{lem3}, which contradicts the maximality of $\rho(G)$. This means that $n_1\geq 2$. Taking a vertex $v_1\in P_1$ and
let $G_1=G-\{uv\mid v\in P_i\setminus \{v_1\},~1\leq i\leq q\}+\{uw_j\mid t+1\leq j\leq\delta-1\}$. According to $x(w_{t+1})>x_1$ and Lemma \ref{lem3}, we get
\begin{equation}\label{equ5}
\begin{split}
\rho(G_1)>\rho(G).
\end{split}
\end{equation}
Furthermore, we take a vertex $v_2\in V(K_{n_1})\setminus \{v_1\}$, let $G_2=G_1+\left\{v_2v \mid v\in V(G)\setminus(S\cup V(K_{n_1}))\right\}$ and $S_1=S-\{u\}+\{v_2\}$. Then we have $|S_1|=\kappa^{\prime}$, $G_2-S_1$ contains $q$ components $(B_1\setminus\{v_2\})\cup\{u\},B_2,\ldots,B_q$ and at least two of them are non-trivial components. Therefore $G_2\in \mathcal{A}_n^{\kappa',\delta}$ and $\rho(G_2)>\rho(G_1)$ by Lemma \ref{lem2}. Combining this with (\ref{equ5}), we have
$$\rho(G_2)>\rho(G),$$
which contradicts the maximality of $\rho(G)$. This indicates that $u\notin S$. \par
\vspace*{2mm}
By Claim \ref{claim1}, $u\notin S$, and so $u\in V(B_i)$ for some $1\leq i\leq q$. Without loss of generality, we may consider that $u\in V(B_1)$.

\begin{claim} $t=\delta-1$.\label{claim2}
\end{claim}
Otherwise, $t=\delta$ or $t<\delta-1$. If $t=\delta$, then $B_1$ is a trivial component, i.e., $V(B_1)=\{u\}$. This means that there must be at least two non-trivial components of $B_2,B_3, \ldots,B_q$ and $q\geq 3$. Without loss of generality, we assume that $B_2$ and $B_3$ are two non-trivial components. By the maximality of $\rho(G)$, we can obtain that $G-u\cong K_{\kappa'}\vee(K_{n_2}\cup K_{n_3}\cup \cdots \cup K_{n_q})$. We can deduce that $q=3$. Otherwise, it is easy to verify that $G-u$ is a proper spanning subgraph of $G^{\star}-u=K_{\kappa'}\vee(K_{n_2}\cup K_{n_3'})$, where $n_3'=\sum_{i=3}^{q}n_i$. According to Lemma \ref{lem2}, we obtain $\rho(G)<\rho(G^{\star})$, which contradicts the maximality of $\rho(G)$. Let $G'$ be the graph obtained from $K_1\cup(K_{\kappa'}\vee(K_{n_2}\cup K_{n_3}))$ by adding $\delta$ edges between the isolated vertex $K_1$ and $K_{\kappa'}$. Hence,
$$\rho(G)\leq \rho(G'),$$
with equality if and only if $G\cong G'$.
Without loss of generality, assume that $n_2\leq n_3$. Notice that $n_2\geq 2$ because $B_2$ is a non-trivial component. Furthermore, let $G''$ be the graph obtained from $K_1\cup(K_{\kappa'}\vee(K_2\cup K_{n-\kappa'-3}))$ by adding $\delta$ edges between the isolated vertex $K_1$ and $K_{\kappa'}$. Combining these with Lemma \ref{lem4}, we have
$$\rho(G')\leq \rho(G''),$$
with equality if and only if $G'\cong G''$. By the maximality of $\rho(G)$, we infer that $G\cong G''$. That is $G-u\cong K_{\kappa'}\vee(K_2\cup K_{n-\kappa'-3})$.  So $(n_2,~n_3)=(2,~n-\kappa'-3)$. Let $N_S(u)=\{w_1,w_2,\ldots,w_\delta\}$, $S\setminus N_S(u)=\{w_{\delta+1},w_{\delta+2},\ldots,w_{\kappa'}\}$ and $V(K_2)=\{z_1,z_2\}$. Let $\bm{x}$ be the Perron vector of $A(G)$ and $\rho=\rho(G)$. By symmetry, $x(v)=x_i$ for $v\in V(K_{n_i})$, where $i=2,~3$. Note that $x(w_1)=x(w_i)$ for $2\leq i\leq \delta$ and $x(w_{\delta+1})=x(w_i)$ for $\delta+2\leq i\leq\kappa'$. Since $N_G(u)\setminus \{w_1\}\subset N_G(w_1)\setminus \{u\}$, $x(w_1)>x(u)$ by Lemma \ref{lem1}. And by $A(G)\bm{x}=\rho \bm{x}$, we obtain
\begin{equation}\label{equ6}
\begin{aligned}
\rho x(u)&=\delta x(w_1),
\end{aligned}
\end{equation}
\begin{equation}\label{equ7}
\begin{aligned}
\rho x_2&=\delta x(w_1)+(\kappa'-\delta)x(w_{\delta+1})+x_2,
\end{aligned}
\end{equation}
\begin{equation}\label{equ8}
\begin{aligned}
\rho x_3&=\delta x(w_1)+(\kappa'-\delta)x(w_{\delta+1})+(n-\kappa'-4)x_3\geq \delta x(w_1)+(n-\kappa'-4)x_3,
\end{aligned}
\end{equation}
\begin{equation}\label{equ9}
\begin{aligned}
\rho x(w_\delta)&=(\delta-1)x(w_1)+(\kappa'-\delta)x(w_{\delta+1})+x(u)+2x_2+(n-\kappa'-3)x_3.
\end{aligned}
\end{equation}
Observe that $G$ contains $K_{n-\kappa'-3}$ as a proper subgraph. Then $\rho>\rho(K_{n-\kappa'-3})=n-\kappa'-4\geq 1$ due to $n-\kappa'-3\geq 2$.
Combining this with (\ref{equ6}) and (\ref{equ8}), we get
\begin{equation}\label{equ10}
\begin{aligned}
x_3\geq \frac{\rho x(u)}{\rho-(n-\kappa'-4)}\geq \frac{\rho x(u)}{\rho-1}.
\end{aligned}
\end{equation}
From (\ref{equ7}), (\ref{equ8}) and (\ref{equ9}), we have
\begin{equation}\label{equ11}
\begin{split}
x_2+x_3-x(w_{\delta})=\frac{\delta x(w_1)+(\kappa'-\delta)x(w_{\delta+1})-x(u)}{\rho+1}.
\end{split}
\end{equation}
Assume that $E_1=\{z_2v|~v\in V(K_{n-\kappa'-3})\}+\{uz_1\}$
and $E_2=\{uw_\delta\}+\{z_1z_2\}$. Let $G_1=G-E_2+E_1$. Clearly, $G_1\cong G_n^{\kappa',\delta}$. Let $\bm{y}$ be the Perron vector of $A(G_1)$ and $\rho'=\rho(G_1)$. By symmetry, $y(v)=y(z_2)$ for $v\in V(K_{n_2'})$, where $n_2'=n-\kappa'-2$. Note that $y(w_1)=y(w_i)$ for $2\leq i\leq \delta-1$ and $y(w_{\delta})=y(w_i)$ for $\delta+1\leq i\leq\kappa'$. Therefore, by $A(G_1)\bm{y}=\rho'\bm{y}$, we have
\begin{equation*}
\begin{split}
\rho'y(z_1)&=(\delta-1)y(w_1)+(\kappa'-\delta+1)y(w_\delta)+y(u),
\end{split}
\end{equation*}
\begin{equation*}
\begin{split}
\rho'y(z_2)&=(\delta-1)y(w_1)+(\kappa'-\delta+1)y(w_\delta)+(n-\kappa'-3)y(z_2),
\end{split}
\end{equation*}
\begin{equation*}
\begin{split}
\rho'y(w_\delta)&=(\delta-1)y(w_1)+(\kappa'-\delta)y(w_\delta)+y(z_1)+(n-\kappa'-2)y(z_2),
\end{split}
\end{equation*}
from which we obtain that
\begin{equation}\label{equ12}
\begin{split}
y(w_\delta)=y(z_2)+\frac{y(z_1)}{\rho'+1}.
\end{split}
\end{equation}
Notice that $N_{G_1}(u)\setminus \{z_1\}\subset N_{G_1}(z_1)\setminus\{u\}$, then $y(z_1)>y(u)$ by Lemma \ref{lem1}. Combining this with $n-\kappa'-4\geq 1$, we get
\begin{align*}
(\rho'-1)(y(z_2)-y(z_1))=(n-\kappa'-4)y(z_2)+y(z_1)-y(u)>0.
\end{align*}
Since $G_1$ contains $K_{n-\kappa'-2}$ as a proper subgraph, it follows that $\rho'>\rho(K_{n-\kappa'-2})=n-\kappa'-3\geq 2$. This suggests that $y(z_2)>y(z_1)$. In view of $y(z_1)>y(u)$, we get $y(z_2)>y(u)$. Combining these with (\ref{equ10}), (\ref{equ11}) and (\ref{equ12}), we have
\begin{equation*}
\begin{split}
&y^{T}(\rho'-\rho)x\\
=&y^{T}(A(G_1)-A(G))x\\
=&\sum_{z_2v\in E_1}(x(z_2)y(v)+x(v)y(z_2))+(x(u)y(z_1)+x(z_1)y(u))-(x(u)y(w_\delta)+x(w_\delta)y(u))\\
&-(x(z_1)y(z_2)+x(z_2)y(z_1))\\
\geq& 2(x_2y(z_2)+x_3y(z_2))+(x(u)y(z_1)+x_2y(u))-(x(u)y(w_\delta)+x(w_\delta)y(u))-(x_2y(z_2)+x_2y(z_1))\\
=&(x_2y(z_2)-x_2y(z_1))+(x_3y(z_2)+x(u)y(z_1)-x(u)y(w_\delta))+(x_3y(z_2)+x_2y(u)-x(w_\delta)y(u))\\
\geq&x_2(y(z_2)-y(z_1))+\left(\frac{\rho x(u)}{\rho-1}y(z_2)+x(u)y(z_1)-x(u)y(w_\delta)\right)+(x_3y(z_2)+x_2y(u)-x(w_\delta)y(u))\\
>&x(u)\left(\frac{\rho}{\rho-1}y(z_2)+y(z_1)-y(w_\delta)\right)+y(u)(x_3+x_2-x(w_\delta))~~(\text{since}~y(z_2)>y(z_1)>y(u))\\
=&x(u)\left(\frac{1}{\rho-1}y(z_2)+\frac{\rho'}{\rho'+1}y(z_1)\right)+y(u)\frac{\delta x(w_1)+(\kappa'-\delta)x(w_{\delta+1})-x(u)}{\rho+1}\\
>&0~~(\text{since}~\rho>1,~\kappa'\geq \delta\geq1~\text{and}~x(w_1)>x(u)).
\end{split}
\end{equation*}
Thus $\rho'>\rho$, which contradicts the maximality of $\rho$. This demonstrates that $t\neq \delta$. Next, we consider $t<\delta-1$. Then $d_{B_1}(u)=\delta-t\geq2$ and $B_1$ is a non-trivial component. Let $N_S(u)=\{w_1,w_2,\ldots,w_t\}$, $S\setminus N_S(u)=\{w_{t+1},w_{t+2},\ldots,w_{\kappa^{\prime}}\}$ and $N_{B_1}(u)=\{v_1,v_2,\ldots,v_{\delta-t}\}$. By the maximality of $\rho(G)$, we can deduce that $G-u\cong K_{\kappa^{\prime}}\vee(K_{n_1-1}\cup K_{n_2}\cup\cdots \cup K_{n_q})$. Assume that $\bm{x}$ is the Perron vector of $A(G)$. By symmetry, $x_1=x(v)$ for $v \in N_{B_1}(u)$, $x_{1}^{\prime}=x(v)$ for $v \in V(B_1) \backslash N_{B_1}[u]$ and $x_i=x(v)$ for $v \in V(K_{n_i})$, where $2\leq i\leq q$. Notice that $x(w_1)=x(w_i)$ for  $2\leq i\leq t$ and $x(w_{t+1})=x(w_i)$ for $t+2 \leq i \leq \kappa^{\prime}$. Then, by $A(G) \bm{x}=\rho(G) \bm{x}$, we have
$$\rho(G)x(w_{t+1})=tx(w_1)+(\kappa^{\prime}-t-1)x(w_{t+1})+(\delta-t)x_1+(n_1-1-(\delta-t)) x_1^{\prime}+\sum_{i=2}^{q}n_ix_i,$$
$$\rho(G)x_1=tx(w_1)+(\kappa^{\prime}-t)x(w_{t+1})+(\delta-t-1)x_1+(n_1-1-(\delta-t))x_1^{\prime}+x(u),$$
from which we obtain that
\begin{equation*}
\begin{split}
&(\rho(G)+1)(x(w_{t+1})-x_1) \\
= &\sum_{i=2}^{q}n_ix_i-x(u) \\
= & \frac{1}{\rho(G)}\left(\sum_{i=2}^{q}n_i(\rho(G)x_i)-\rho(G)x(u)\right) \\
= & \frac{\sum_{i=2}^{q}n_i(tx(w_1)+(\kappa'-t)x(w_{t+1})+(n_i-1)x_i)-(tx(w_1)+(\delta-t)x_1)}{\rho(G)} \\
> & \frac{(\sum_{i=2}^{q}n_i-1)tx(w_1)+\sum_{i=2}^{q}n_i(\kappa'-t)x(w_{t+1})-(\delta-t)x_1}{\rho(G)}~~(\text{since}~\sum_{i=2}^{q}n_i\geq q) \\
\geq & \frac{(q-1)tx(w_1)+q(\kappa'-t)x(w_{t+1})-(\delta-t)x_1}{\rho(G)}~~(\text{since}~\sum_{i=2}^{q}n_i\geq q) \\
>&\frac{(\delta-t)(x(w_{t+1})-x_1)}{\rho(G)}~~(\text{since}~ q\geq 2~\text{and}~\kappa^{\prime}\geq\delta>t\geq 0).
\end{split}
\end{equation*}
It follows that $\frac{(\rho^{2}(G)+\rho(G)-\delta+t)(x(w_{t+1})-x_1)}{\rho(G)}>0$. Note that $G$ contains $K_{\kappa^{\prime}+n_2}$ as a proper subgraph. Thus, $\rho(G)>\rho(K_{\kappa^{\prime}+n_2})=\kappa^{\prime}+n_2-1\geq\delta$ due to $\kappa^{\prime}\geq\delta$ and $n_2\geq1$. Combining this with $t\geq 0$, we can deduce that $\rho^{2}(G)+\rho(G)-\delta+t>0$. This suggests that $x(w_{t+1})>x_1$. Let $G_2=G-\{uv_i|~2\leq i \leq\delta-t\}+\{uw_j|t+1\leq j\leq\delta-1\}$. It is obvious that $G_2\in \mathcal{A}_n^{\kappa',\delta}$. According to Lemma \ref{lem3}, we have $\rho(G_2)>\rho(G)$, which contradicts the maximality of $\rho(G)$. Therefore, we have $t=\delta-1$, completing the proof of Claim \ref{claim2}.

\vspace*{2mm}
Without loss of generality, we still assume that $u\in V(B_1)$ and $d_{B_1}(u)=1$ by Claim \ref{claim2}. For convenience, we denote $z=N_{B_1}(u)$. Then $B_1$ is a non-trivial component. By the maximality of $\rho(G)$, we can deduce that $G-u\cong K_{\kappa'}\vee(K_{n_1-1}\cup K_{n_2}\cup \cdots \cup K_{n_q})$. Furthermore, we can conclude that $q=2$. Otherwise, let $G^{*}=G+\{vv'|~v\in V(B_2),~v'\in V(B_i),~3\leq i\leq q\}$. One can easily observe that $G^{*}\in \mathcal{A}_n^{\kappa',\delta}$ and $G$ is a proper subgraph of $G^{*}$. According to Lemma \ref{lem2}, we have $\rho(G)<\rho(G^{*})$, which contradicts the maximality of $\rho(G)$. Therefore, $G-S$ has exactly two non-trivial components $B_1$ and $B_2$. Recall that $|V(B_i)|=n_i$, then $n_i\geq 2$ for $i=1,~2$.

\begin{claim} $n_1\leq n_2$.\label{claim3}
\end{claim}
Otherwise, $n_1>n_2$. We can deduce that $N_{B_1}(z)\setminus \{u\}\neq \emptyset$, that is $d_{B_1}(z)\geq2$. Let $N_S(u)=\{w_1,w_2,\ldots,w_{\delta-1}\}$, $S\setminus N_S(u)=\{w_{\delta},w_{\delta+1},\ldots,w_{\kappa^{\prime}}\}$ and $V(B_2)=\{u_1,u_2,\ldots,u_{n_2}\}$. By the maximality of $\rho(G)$, we infer that $G-u\cong K_{\kappa^{\prime}}\vee(K_{n_1-1}\cup K_{n_2})$. Assume that $\bm{x}$ is the Perron vector of $A(G)$. By symmetry, $x(v)=x_1$ for $v\in N_{B_1}(z)\setminus\{u\}$ and $x(u_i)=x_2$ for $1\leq i\leq n_2$. Note that $x(w_1)=x(w_i)$ for $2\leq i\leq \delta-1$ and $x(w_\delta)=x(w_i)$ for $\delta+1\leq i\leq \kappa^{\prime}$. Therefore, by $A(G)\bm{x}=\rho(G)\bm{x}$, we have
\begin{equation}\label{equ13}
\begin{split}
\rho(G)x(z)&=(\delta-1)x(w_1)+(\kappa^{\prime}-\delta+1)x(w_\delta)+(n_1-2)x_1+x(u),
\end{split}
\end{equation}
\begin{equation}\label{equ14}
\begin{split}
\rho(G)x_1&=(\delta-1)x(w_1)+(\kappa^{\prime}-\delta+1)x(w_\delta)+(n_1-3)x_1+x(z),
\end{split}
\end{equation}
\begin{equation}\label{equ15}
\begin{split}
\rho(G)x_2&=(\delta-1)x(w_1)+(\kappa^{\prime}-\delta+1)x(w_\delta)+(n_2-1)x_2,
\end{split}
\end{equation}
\begin{equation}\label{equ16}
\begin{split}
\rho(G)x(w_\delta)&=(\delta-1)x(w_1)+(\kappa'-\delta)x(w_\delta)+(n_1-2)x_1+x(z)+n_2x_2.
\end{split}
\end{equation}
Notice that $N_G(v)\setminus \{z\}\subset N_G(z)\setminus\{v\}$ for $v\in N_{B_1}(z)\setminus\{u\}$, we have $x(z)>x_1$ by Lemma \ref{lem1}. Similarly, $x(z)>x(u)$. From (\ref{equ14}) and (\ref{equ15}), we get
$$(\rho(G)-n_2+1)(x_1-x_2)=(n_1-n_2-2)x_1+x(z)>0$$
due to $n_1\geq n_2+1$ and $x(z)>x_1$. Since $G$ contains $K_{\kappa^{\prime}+n_2}$ as a proper subgraph and $\kappa'\geq 1$, it follows that $\rho(G)>\rho(K_{\kappa^{\prime}+n_2})=\kappa^{\prime}+n_2-1\geq n_2$, and hence $x_1>x_2$.
Suppose that $E_1'=\{u_{n_2}v\mid v\in V(B_1)\setminus \{u\}\}+\{uu_1\}$ and $E_2'=\{u_{n_2}u_i\mid 1\leq i\leq n_2-1\}+\{zu\}$. Let $G_1=G+E_1'-E_2'$. Then $G_1\in \mathcal{A}_n^{\kappa',\delta}$ and we have
\begin{equation*}
\begin{aligned}
&\rho(G_1)-\rho(G)\\
\geq& \bm{x}^{T}(A(G_1)-A(G))\bm{x}\\
=&2\sum_{u_{n_2}v\in E_1'}x(u_{n_2})x(v)+2x(u)x(u_1)-2\sum_{u_{n_2}u_i\in E_2'}x(u_{n_2})x(u_i)-2x(z)x(u)\\
=&2(n_1-2)x_1x_2+2x(z)x_2+2x(u)x_2-2(n_2-1)x_2x_2-2x(z)x(u)\\
>&2x(z)x_2-2x(z)x(u)~~(\text{since}~n_1\geq n_2+1~\text{and}~x_1>x_2)\\
=&\frac{2x(z)}{\rho(G)}(\rho(G)x_2-\rho(G)x(u))\\
=&\frac{2x(z)[(\delta-1)x(w_1)+(\kappa'-\delta+1)x(w_\delta)+(n_2-1)x_2-(\delta-1)x(w_1)-x(z)]}{\rho(G)}\\
\geq&\frac{2x(z)(x(w_\delta)+x_2-x(z))}{\rho(G)}~~(\text{since}~\kappa'\geq\delta~\text{and}~n_2\geq 2)\\
=&\frac{2x(z)(\rho(G)x(w_\delta)+\rho(G)x_2-\rho(G)x(z))}{\rho^2(G)}\\
=&\frac{2x(z)[(\delta-1)x(w_1)+(\kappa'-\delta)x(w_\delta)+(2n_2-1)x_2+x(z)-x(u)]}{\rho^2(G)}~~(\text{by}~(\ref{equ13}),~(\ref{equ15})~\text{and}~(\ref{equ16}))\\
>&0~~(\text{since}~\delta\geq 1,~\kappa'\geq\delta,~n_2\geq 2~\text{and}~x(z)>x(u)).
\end{aligned}
\end{equation*}
Hence $\rho(G_1)>\rho(G)$, which contradicts the maximality of $\rho(G)$. This implies that $n_1\leq n_2$, proving Claim \ref{claim3}.

\vspace*{2mm}
In what follows, we shall prove that $G\cong G_n^{\kappa',\delta}$ for $\kappa'>\delta-1$. In fact, by Claims \ref{claim1}--\ref{claim3} and the maximality of $\rho(G)$, we can infer that $G-u\cong K_{\kappa'}\vee(K_{n_1-1}\cup K_{n_2})$ for $2\leq n_1\leq n_2$.
If $n_1=2$ and $n_2=n-\kappa'-2$, then the result follows. If not, by (\ref{equ13}) and (\ref{equ15}), we get
$$(\rho(G)-n_2+1)(x_2-x(z))=(n_2-1)x(z)-(n_1-2)x_1-x(u)>0$$
due to $n_2\geq n_1\geq3$, $x(z)>x_1$ and $x(z)>x(u)$. Recall that $\rho(G)>n_2$, then $x_2>x(z)$.
Assume that $E_1''=\{vv'\mid v\in N_{B_1}(z)\setminus\{u\}, v'\in V(K_{n_2})\}$ and $E_2''=\{vz\mid v\in N_{B_1}(z)\setminus\{u\}\}$. Let $G_2=G-E_2''+E_1''$. Clearly, $G_2\in \mathcal{A}_n^{\kappa',\delta}$ and
\begin{equation*}
\begin{aligned}
\rho(G_2)-\rho(G) \geq& \bm{x}^{T}(A(G_1)-A(G))\bm{x}\\
=& 2\sum_{vv'\in E_1''}x(v)x(v')-2\sum_{vz\in E_2''}x(v)x(z)\\
=& 2(n_1-2)n_2x_1x_2-2(n_1-2)x_1x(z)\\
=& 2(n_1-2)x_1(n_2x_2-x(z))\\
>&0~~(\text{since}~n_2\geq n_1\geq3~\text{and}~x_2>x(z)).
\end{aligned}
\end{equation*}
It follows that $\rho(G_2)>\rho(G)$, which contradicts the maximality of $\rho(G)$. This means that $n_1=2$.
Then we conclude that $G\cong G_n^{\kappa^{\prime}, \delta}$, as desired.

\vspace*{2mm}
\begin{case} $\kappa^{\prime}\leq\delta-1$.
\end{case}
We first assert that $q=2$. Otherwise, $q\geq 3$. Without loss of generality, assume that $B_1$ is one of the non-trivial components. Let $G^{\prime}=K_{\kappa'}\vee(K_{n_1}\cup K_{n_2'})$, where $n_2'=\sum_{i=2}^{q}n_i$. It is not difficult to see that $G$ is a proper spanning subgraph of $G^{\prime}$. Thus $\rho(G)<\rho(G^{\prime})$ by Lemma \ref{lem2}, which contradicts the maximality of $\rho(G)$. This implies that $q=2$. Therefore, $G-S$ has exactly two non-trivial components $B_1$ and $B_2$. Furthermore,
$$\rho(G)\leq \rho(K_{\kappa'}\vee(K_{n_1}\cup K_{n_2})),$$
with equality if and only if $G\cong K_{\kappa'}\vee(K_{n_1}\cup K_{n_2})$. Without loss of generality,
we may assume that $n_1\geq n_2$. Evidently, $n_2\geq \delta-\kappa^{\prime}+1$ because the minimum degree of $G$ is $\delta$. Combining this with Lemma \ref{lem4}, we get
$$\rho(K_{\kappa'}\vee(K_{n_1}\cup K_{n_2}))\leq \rho(K_{\kappa'}\vee(K_{n-\delta-1}\cup K_{\delta-\kappa'+1})),$$
with equality if and only if $(n_1,n_2)=(n-\delta-1,\delta-\kappa'+1)$. By the maximality of $\rho(G)$, we infer that $G\cong K_{\kappa'}\vee(K_{n-\delta-1}\cup K_{\delta-\kappa'+1})$.\par
This completes the proof.

$\hfill\qedsymbol$

\section{Proof of Theorem \ref{thm2}} \label{3}
Let $M$ be a real matrix of order $n$, and give a partition $\pi=(V_1,V_2,\ldots,V_k)$ of $V(G)$, where $|V_i|=n_i$. The matrix $M$ described in the following
block form
$$M=\left(
             \begin{array}{ccccc}
               M_{11} & M_{12} & \cdots & M_{1k} \\
               M_{21} & M_{22} & \cdots & M_{2k}\\
               \vdots & \vdots & \vdots & \vdots\\
               M_{k1} & M_{k2} & \cdots & M_{kk}
             \end{array}
           \right),$$
where the blocks $M_{ij}$ are $n_i\times n_j$ matrices for any $1\leq i,j \leq k$ and $n = n_1+n_2+\cdots+n_k$.
For $1\leq i,j \leq k$, let $b_{ij}$ denote the average value of all row sums of $M_{ij}$, i.e., $b_{ij}$ is the sum of all
entries in $M_{ij}$ divided by the number of rows. Then $B_\pi(M)= (b_{ij})_{i,j=1}^{k}$ (or simply $B_\pi$) is
called the quotient matrix of $M$ with respect to $\pi$. Furthermore, if each block $M_{ij}$ of $M$ has a constant row
sum, then $B_\pi$ is called the equitable quotient matrix of $M$.

\begin{lemma}\label{lem5} \cite{L.H. You}
Let $M$ be a real symmetric matrix, and let $\lambda_1(M)$ be the largest eigenvalue of $M$. If $B_\pi$ is an equitable quotient matrix of $M$, then the eigenvalues of $B_\pi$ are also eigenvalues of $M$. Furthermore, if $M$ is a nonnegative irreducible matrix, then $\lambda_1(M)=\lambda_1(B_\pi)$.
\end{lemma}

\begin{lemma} \cite{J.A. Bondy}\label{lem6}
Let $D$ be an arbitrary strongly connected digraph with vertex connectivity $k$. Suppose that $S$ is a $k$-vertex cut of $D$ and $D_1,D_2,\ldots,D_s$ are the strongly connected components of $D-S$. Then there exists an ordering of $D_1,D_2,\ldots,D_s$ such that, for $1\leq i\leq s$ and $v\in V(D_i)$, every tail of $v$ in $D_1,D_2,\ldots,D_{i-1}$.
\end{lemma}

\begin{lemma}\label{lem7} \cite{W.X. Hong}
Let $D$ be a strongly connected digraph and $D'$ be a proper subgraph of $D$. Then $\rho(D')<\rho(D)$.
\end{lemma}

\begin{lemma}\label{lem8}
Let $f(x)=4x^2-4(n-k)x+n^2$, where $2\leq x\leq n-k-2$ and $k\geq 1$. Then $f_{max}(x)=f(2)=f(n-k-2)=(n-4)^2+8k>0$.
\end{lemma}

\vspace*{2mm}
\noindent
{\bf{Proof of Theorem \ref{thm2}}.} Suppose that $D$ is a digraph with maximum spectral radius in $\mathcal{D}_{n, k}$. Thus, there exists some subset $S\subseteq V(D)$ with $|S|=k$ such that $D-S$ contains at least two strong connected non-trivial components. By Lemma \ref{lem6},
we know that $D_1$ with $|V(D_1)|=m$ is the strongly connected non-trivial component of $D-S$ where $N^-(v)=\{u\in V(D-S-D_1)~|~(u,v)\in E(D)\}$ and $d^-(v)=|N^-(v)|=0$ for any $v\in V(D_1)$. Let $D_2 = D-S-D_1$. Since $D$ is a strongly connected digraph, we add arcs to $D$ until both induced subdigraph of $V(D_1)\cup S$ and induced subdigraph of $V(D_2)\cup S$ attain to complete digraphs, add arc
$(u,v)$ for any $u\in V(D_1)$ and any $v\in V (D_2)$. Denote the new digraph by $D'$. It is easy to find that $D'=\vec{G}_n^{k,m}\in \vec{\mathcal{G}}(n,k)\subseteq \mathcal{D}_{n,k}$.
And $D$ is the subdigraph of $D'$, then $\rho(D')\geq\rho(D)$, with equality if and only if $D\cong D'$ by Lemma \ref{lem7}. Therefore, the digraphs which achieve the maximum spectral radii among all digraphs in $\mathcal{D}_{n,k}$ must be in $\vec{\mathcal{G}}(n,k)$.
By the maximality of $\rho(D)$, we have $D\cong \vec G_{n}^{k,m}$ for an integer $m\geq 2$ and $n-k-m\geq 2$. Observe that $A(D)$ has the equitable quotient matrix
$$B_\pi=\left(
             \begin{array}{ccccc}
               k-1 & m & n-k-m \\
               k & m-1 & n-k-m\\
               k & 0 & n-k-m-1
             \end{array}
           \right).$$
Then
\begin{align*}
\varphi(B_\pi,\lambda)&=\begin{vmatrix}
                        \lambda-(k-1) & -m & -(n-k-m) \\
                        -k & \lambda-(m-1) & -(n-k-m) \\
                        -k &  0 &  \lambda-(n-k-m-1)
                        \end{vmatrix}\\
\ &=\begin{vmatrix}
                       \lambda+1 & -(\lambda+1) & 0 \\
                        0 & \lambda-(m-1) & -(\lambda+1) \\
                        -k & 0 & \lambda-(n-k-m-1)
                        \end{vmatrix} \\
\ &=(\lambda+1)[(\lambda-(m-1))(\lambda-(n-k-m-1))-k(\lambda+1)]\\
\ &=(\lambda+1)[\lambda^2+(2-n)\lambda+mn-mk-m^2-n+1].
\end{align*}
Combining this with Lemma \ref{lem5}, it shows that $\rho(D)$ is equal to the largest root of the quadratic equation $x^2+(2-n)x+mn-mk-m^2-n+1=0$, where $2\leq m\leq n-k-2$. Thus we have $\rho(D)=\frac{n-2+\sqrt{4m^2-4(n-k)m+n^2}}{2}\leq \frac{n-2+\sqrt{(n-4)^2+8k}}{2}$ by Lemma \ref{lem8}, equality holding if and only if $m=2$ or $m=n-k-2$. If $m=2$, then $D\cong \vec {G}_n^{k,2}$; if $m=n-k-2$, then $D\cong \vec {G}_n^{k,n-k-2}$, as required.\par
This completes the proof.

$\hfill\qedsymbol$

\end{document}